\documentclass[letterpaper,10pt,conference]{ieeeconf}

\IEEEoverridecommandlockouts
\usepackage{amsmath,amssymb,amsfonts}
\usepackage{xurl}
\usepackage[colorlinks=true,urlcolor=blue,linkcolor=black,citecolor=black]{hyperref}
\usepackage{graphicx}
\usepackage{booktabs}
\usepackage{makecell}
\usepackage{cite}
\usepackage[caption=false,font=footnotesize]{subfig}
\usepackage{algorithm}
\usepackage{algpseudocode}
\usepackage{stfloats}

\title{\LARGE \bf
Backstepping-Guided Reinforcement Learning for Wide-Range Saint-Venant Canal Regulation
}

\author{
Chenchen Wang$^{1}$ and Jie Qi$^{2}$
\thanks{The paper is supported by the National Natural Science Foundation of China
(62173084).}
\thanks{$^{1}$ Chenchen Wang and $^{2}$Jie Qi are with the School of Information
and Intelligent Science, Donghua University, Shanghai, China
(e-mail: 2242131@mail.dhu.edu.cn; jieqi@dhu.edu.cn).
Corresponding author: Jie Qi.}%
}
\begin{document}

\maketitle
\thispagestyle{empty}
\pagestyle{empty}

\begin{abstract}
Backstepping control provides local stability guarantees for nonlinear Saint-Venant systems,
but its regulation performance may degrade when the system operates far from the nominal equilibrium. This letter proposes a backstepping-guided soft
actor-critic (SAC) controller framework that incorporates model-based control knowledge into reinforcement learning (RL). The nominal backstepping control law is first learned by deep operator network (DeepONet) and embedded into the actor and critic networks as prior-informed feature representations. The learned prior is further combined with the SAC policy to generate the final control input, while a transfer-learning strategy preserves the useful backstepping knowledge during adaptation to the nonlinear dynamics.
Simulation results on the Sambre River model demonstrate that the
proposed method improves learning efficiency and maintains effective regulation over larger initial deviations than backstepping control. 
% Nonlinear Saint-Venant canal-flow systems are difficult to regulate due to distributed hyperbolic dynamics, boundary actuation, nonlinear friction, and long-distance transport. Nominal backstepping provides a stabilizing prior for the linearized model, but its transient performance may degrade under large off-nominal initial conditions. This letter proposes a backstepping-guided residual deep operator network-soft actor-critic (DeepONet-SAC) controller. The nominal backstepping gate command is encoded into an actor-side DeepONet to provide a baseline command and policy features. SAC learns a correction guided by this baseline, while an independent critic-side encoder from the same pretrained DeepONet provides twin-Q state features. A partial parameter-freezing transfer strategy is adopted, while the remaining parameters are fine-tuned with a smaller learning rate.
% Sambre River model simulations show that, under near-equilibrium initial conditions, the proposed method maintains good control performance; under large initial deviation conditions, it improves training behavior and reduces state-error responses.

\end{abstract}


\section{Introduction}
\label{sec:introduction}

Open-channel flow regulation plays a crucial role in hydraulic
engineering, irrigation, and water-resource management. The
Saint-Venant equations are widely used to describe the coupled dynamics of water depth and discharge through a nonlinear hyperbolic partial differential equation (PDE) model \cite{dosSantos2008boundary,hayat2022pi}. However, the nonlinear and spatially distributed nature of canal systems makes boundary control over a wide range of operating conditions particularly challenging. 

Backstepping provides a systematic model-based framework for boundary
control of Saint-Venant systems. By constructing integral
transformations, backstepping controllers can be designed based on
linearized models and guarantee local stability of the closed-loop
nonlinear systems around prescribed operating equilibria
\cite{dosSantos2008boundary,hayat2022pi,coron2013local}. However, their regulation performance may deteriorate when the nonlinear system is initialized far from the equilibrium. This motivates the development of control strategies that extend effective regulation to larger initial deviations.

Reinforcement learning (RL) can optimize control policies directly from nonlinear system trajectories and thus offers a potential means of improving regulation beyond the neighborhood of nominal operating conditions. An early study addressed PDE control with high-dimensional continuous action spaces by exploiting spatial regularities among action dimensions \cite{pan2018reinforcement}. In fluid systems governed by nonlinear PDEs, 
 % deep RL has been successfully applied to active flow control\cite{rabault2019artificial}, while 
stability and sensitivity information has been incorporated into RL-based control of cylinder wakes \cite{li2022reinforcement}. For hyperbolic PDEs, RL-based boundary control has been investigated for congested traffic-flow models and compared with PDE backstepping and PI control \cite{yu2022rl}. More recently, benchmark environments such as Controlgym and PDE Control Gym have facilitated systematic evaluation of RL algorithms on infinite-dimensional PDE control problems \cite{zhang2024controlgym,bhan2024pde}. Convolutional architectures have been introduced to process high-dimensional distributed states \cite{peitz2024distributed}, while online RL methods have been investigated for PDE identification and control under unknown dynamics \cite{alla2024online}. 
Surrogate-model-assisted RL has further been shown to reduce the amount of training data required from the original PDE system \cite{werner2023learning}. 
Nevertheless, model-free RL for PDE boundary control generally requires extensive trajectory generation and makes limited use of available analytical control knowledge. Recent physics-informed actor--critic approaches have shown that incorporating Lyapunov-based control structures into policy learning can improve stability-oriented controller design \cite{wang2024actorcritic}.

Therefore, this letter incorporates backstepping control knowledge into a soft actor-critic (SAC) framework \cite{haarnoja2018soft} to improve training efficiency and nonlinear canal-flow regulation. A neural operator, DeepONet \cite{lu2021learning1}, is employed to
approximate the backstepping feedback operator, motivated by recent advances in neural-operator approximation of PDE backstepping kernels and control laws 
\cite{bhan2024neural,krstic2024neural,qi2024neural,11180047}. The
pretrained DeepONet provides a backstepping-informed stabilizing prior,
while SAC learns nonlinear corrections to improve regulation under
large initial deviations. 

The main contributions are summarized as follows:
1) The backstepping
feedback operator is approximated by DeepONet and then embedded
into both the actor and critic networks of SAC as prior-informed feature
extractors. A transfer learning strategy for the DeepONet with partially frozen
layers and reduced learning rates is adopted to retain the backstepping prior and improve learning efficiency.
2)  The SAC policy is combined with the learned backstepping prior to compensate for nonlinear effects, which extends the effective regulation range to larger initial deviations than those handled by the backstepping controller.

The remainder of this letter is organized as follows. Section II
presents the Saint-Venant model and the backstepping prior. Section III
introduces the proposed backstepping-guided DeepONet-SAC controller.
Section IV presents numerical results, and Section V concludes the
letter.

\section{Preliminaries}
\label{sec:preliminaries}

\subsection{Saint-Venant Canal Model}

Consider an open-channel canal on $x\in[0,L]$, governed by the Saint-Venant equations~\cite{dosSantos2008boundary}:
\begin{align}
\label{eq:fx}
\begin{aligned}
&\partial_t H+\frac{1}{B}\partial_x Q=0,\\
&\partial_t Q
+\frac{1}{B}\partial_x\left(\frac{Q^2}{H}\right)
+\frac{gB}{2}\partial_x H^2
=gBH\bigl(I-J(H,Q)\bigr),
\end{aligned}
\end{align}
where $H$ and $Q$ are state variables denoting the water depth and discharge, respectively. $B$ and $g$ are constants representing the channel width and gravitational acceleration, and $I$ denotes the bottom slope, which may be constant or spatially varying. The friction slope is
\begin{equation}\label{eq:J}
J(H,Q)=\frac{Q^2}{K^2(BH)^2R^{4/3}(H)},\qquad
R(H)=\frac{BH}{B+2H},
\end{equation}
with $K$ the Strickler coefficient.

The upstream gate is fixed at $U_0(t)=\bar U_0$, while the downstream gate provides the control input:
\begin{align}
Q(t,0)
&=\bar U_0B\mu_0
\sqrt{2g\bigl(z_{\rm up}-H(t,0)\bigr)},
\label{eq:Q_bnd}\\
H(t,L)
&=
\left(
\frac{Q^2(t,L)}{2gB^2\mu_L^2}
\right)^{1/3}
+h_s+U_L(t).
\label{eq:H_bnd}
\end{align}

For a desired steady state $(\bar H(x),\bar Q)$, define
\begin{equation}
\label{eq:wucha}
e_H=H-\bar H,\qquad e_Q=Q-\bar Q.
\end{equation}
The objective is to regulate $(e_H,e_Q)$ to zero through $U_L(t)$.

\subsection{Backstepping Controller}

Under a subcritical steady operating condition, the linearized system can be transformed into characteristic coordinates
\begin{equation}
\label{eq:w}
w(t,x)=
\begin{bmatrix}u(t,x)\\
v(t,x)\end{bmatrix}
=
T(x)
\begin{bmatrix}e_H(t,x)\\
e_Q(t,x)\end{bmatrix},
\end{equation}
where
\begin{equation}
\label{eq:T_matrix}
T(x)=\bigl[t_{ij}(x)\bigr]_{i,j=1}^{2}
=
\begin{bmatrix}
\dfrac{\bar c(x)-\bar V(x)}{2\bar c(x)}
&
\dfrac{1}{2B\bar c(x)}
\\[2mm]
\dfrac{\bar c(x)+\bar V(x)}{2\bar c(x)}
&
-\dfrac{1}{2B\bar c(x)}
\end{bmatrix},
\end{equation}
with
\begin{equation}
\label{eq:Vc}
\bar V(x)=\frac{\bar Q}{B\bar H(x)},\qquad
\bar c(x)=\sqrt{g\bar H(x)}.
\end{equation}

Following~\cite{coron2013local}, the backstepping boundary feedback is
\begin{equation}
\label{eq:bs_characteristic}
U_{\rm err}(t)
=
\int_0^L K_{vu}(L,y)u(t,y)dy
+
\int_0^L K_{vv}(L,y)v(t,y)dy,
\end{equation}
where $K_{vu}$ and $K_{vv}$ denote the backstepping kernels given in~\cite{coron2013local}.

Using the downstream boundary condition, the corresponding gate command is
\begin{equation}
\label{eq:bs_gate}
\begin{aligned}
U_{\rm bs}(t)
=&\bar H(L)
+\frac{
U_{\rm err}(t)
-t_{22}(L)\bigl(Q(t,L)-\bar Q\bigr)}
{t_{21}(L)}
\\
&-h_s
-\left(
\frac{Q^2(t,L)}
{2gB^2\mu_L^2}
\right)^{1/3}.
\end{aligned}
\end{equation}
This controller locally stabilizes the nonlinear Saint-Venant system around $(\bar H,\bar Q)$ under the assumptions of~\cite{coron2013local} and serves as the control prior for the learning-based controller.

\section{DeepONet-SAC control framework}
\begin{figure}[htbp]
	\centering
	\includegraphics[width=1\linewidth]{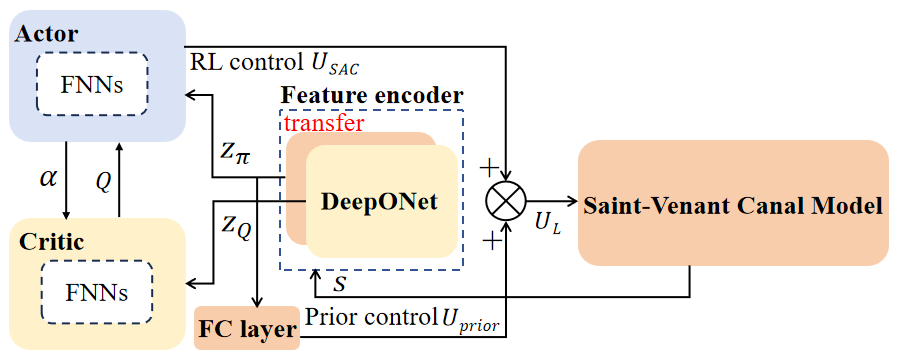}
	\caption{The framework of the transfer DeepONet Prior SAC (TDPSAC). 
    % Actor DeepONet outputs baseline/policy features; critic DeepONet outputs twin-Q features.
    }
	\label{fig:deeponet_sac}
\end{figure}

 %------------------------------------------------
We introduce a DeepONet-SAC control framework, termed Transferred DeepONet Prior SAC (TDPSAC), as illustrated in Fig.~\ref{fig:deeponet_sac}. The proposed framework
transfers backstepping control knowledge into SAC through three
pretrained DeepONet instances, which replace the conventional feature
extractors of the actor, critic, and target critic networks. The
DeepONets are pretrained to approximate the backstepping feedback
operator mapping the distributed Saint-Venant state profile to the
boundary control input.

The transferred DeepONet modules provide prior-informed representations
for SAC and are further adapted through a prior-preserving transfer
strategy. Selected DeepONet layers are frozen, while the remaining
parameters are fine-tuned with reduced learning rates to retain
backstepping knowledge during adaptation to nonlinear dynamics. Moreover, the actor-side DeepONet is followed by an output head
consisting of three fully connected layers to generate the prior control
signal $U_{\mathrm{prior}}(t)$.
The final boundary input is constructed as
\begin{equation}\label{eq:final_U}
U_L(t)=\rho(t)U_{\mathrm{SAC}}(t)+(1-\rho(t))U_{\mathrm{prior}}(t),
\end{equation}
where $\rho(t)\in[0,1]$ is a blending factor balancing the SAC policy and
the backstepping prior. Through knowledge transfer at both the feature
representation and control levels, TDPSAC combines the model-based
backstepping prior with the adaptability of reinforcement learning.

%-----------------------------------------------------new subsection
\subsection{Backstepping-Pretrained DeepONet}

The analytical backstepping controller provides a locally stabilizing
feedback law for the Saint-Venant system around a prescribed equilibrium.
To transfer this model-based control knowledge to the reinforcement learning
framework, a DeepONet is first pretrained to approximate the backstepping
state-feedback mapping.

At each control step, the distributed canal states are sampled at
$x_i\in[0,L]$, $i=1,\ldots,N$, and collected as
\begin{equation}
s=
\left[
H(x_1),\ldots,H(x_N),
Q(x_1),\ldots,Q(x_N)
\right]^{\top}.
\label{eq:deeponet_state}
\end{equation}
The analytical
backstepping controller defines the mapping
\begin{equation}
U_{\mathrm{bs}}
=
\mathcal{K}_{\mathrm{bs}}(s),
\label{eq:bs_operator}
\end{equation}
which is approximated by a DeepONet
\begin{equation}
\mathcal{G}_{\theta}:s\mapsto U_{\mathcal{G},\mathrm{bs}}.
\label{eq:deeponet_operator}
\end{equation}

The pretraining dataset is generated from closed-loop trajectories of the
Saint-Venant system under the analytical backstepping controller \eqref{eq:bs_gate}. The initial conditions are randomly perturbed around the prescribed
equilibrium using Chebyshev functions \cite{qi2024neural}.
% :
% \begin{equation}
% \begin{aligned}
% H_{\ell}^{0}(x)
% &=
% \bar H(x)
% +
% \sum_{j=0}^{M}
% a_{\ell,j}
% T_j\!\left(\frac{2x}{L}-1\right),\\
% Q_{\ell}^{0}(x)
% &=
% \bar Q(x)
% +
% \sum_{j=0}^{M}
% b_{\ell,j}
% T_j\!\left(\frac{2x}{L}-1\right),
% \end{aligned}
% \label{eq:random_initial_profiles}
% \end{equation}
% where $T_j$ denotes the $j$th Chebyshev polynomial. The coefficients
% $a_{\ell,j}$ and $b_{\ell,j}$ are independently sampled from zero-mean
% Gaussian distributions whose standard deviations decay with the mode
% index $j$ as $(j+1)^{-2}$. The resulting Chebyshev expansions define
% the initial-condition perturbations relative to the steady reference
% profiles, which are subsequently normalized by their infinity norms
% and rescaled to the prescribed amplitude ranges.
% , with their ranges chosen such that the resulting initial
% profiles remain within the admissible operating region.

For each sampled initial condition, the closed-loop system is simulated
under the backstepping controller to generate a trajectory.
% \begin{equation}
% \left\{
% \left(s_{\ell}^{k},U_{\mathrm{bs},\ell}^{k}\right)
% \right\}_{k=0}^{K_{\ell}},
% \end{equation}
% where $k$ denote the step index for the trajectory.  
The state-control pairs collected from all trajectories are pooled, randomly
shuffled, and reindexed to form the supervised pretraining dataset
\begin{equation}
\mathcal{D}_{\mathrm{bs}}
=
\left\{
\left(s^{(n)},U_{\mathrm{bs}}^{(n)}\right)
\right\}_{n=1}^{N_{\mathrm{bs}}}.
\label{eq:bs_dataset}
\end{equation}
The DeepONet parameters are obtained by minimizing
\begin{equation}
\theta_D
=
\arg\min_{\theta}
\frac{1}{N_{\mathrm{bs}}}
\sum_{n=1}^{N_{\mathrm{bs}}}
\left|
\mathcal{G}_{\theta}\!\left(s^{(n)}\right)
-
U_{\mathrm{bs}}^{(n)}
\right|^2.
\label{eq:deeponet_pretraining_loss}
\end{equation}

\subsection{Transfer DeepONet Prior SAC}
In particular, the subscripts $\pi$ and $Q$ denote the actor and critic
networks, respectively. The DeepONet-based latent features are defined as
\begin{equation}
z_{\pi}=\Phi_{\pi}(s;\theta_{\pi}),\qquad
z_Q=\Phi_Q(s;\theta_Q),
\label{eq:latent_features}
\end{equation}
where $\Phi_{\pi}$ and $\Phi_Q$ denote the corresponding DeepONet mappings
from the system state to latent feature vectors. The actor-side DeepONet
retains its pretrained control head. Therefore, in addition to providing
$z_{\pi}$ to the SAC actor, it generates the prior control signal
\begin{equation}
U_{\mathrm{prior}}
=
\mathcal{G}_{\theta_{\pi}}(s)
=
W_{\rm p}z_\pi+b_{\rm p},
\label{eq:prior_control}
\end{equation}
where $\mathcal{G}_{\theta_{\pi}}(s)$, defined in
\eqref{eq:deeponet_operator}, denotes the pretrained DeepONet control
mapping that approximates the backstepping feedback operator.

% \begin{figure}[htbp]
% 	\centering
% 	\includegraphics[width=1\linewidth]{pic/312.png}
% 	\caption{Actor-side DeepONet architecture. 
%     % Actor DeepONet outputs baseline/policy features; critic DeepONet outputs twin-Q features.
%     }
% 	\label{fig:deeponet_action}
% \end{figure}

The SAC actor generates the stochastic control component according to
\begin{equation}
U_{\rm SAC}\sim\pi_{\phi}(\cdot|z_\pi),
\end{equation}
and the boundary input applied to the plant is constructed by blending
the SAC output with the backstepping prior according to
\eqref{eq:final_U}. 
% Therefore, the critic evaluates the composite action $U_L$ rather than $U_{\rm SAC}$ alone.

A transfer-learning strategy is further introduced for the DeepONet
parameters. For the actor and online critic DeepONet networks, the
parameters are partitioned as
\begin{equation}
\theta_D=(\theta_{f,D},\theta_{s,D}),
\label{eq:parameter_partition}
\end{equation}
where $\theta_{f,D}$ is frozen during SAC training and $\theta_{s,D}$ is
fine-tuned with a reduced learning rate $\eta_s$. Specifically, the frozen
part contains the fully connected layers except the last layer of the trunk
network and the convolutional layers of the branch network, while the
remaining layers are allowed to adapt. This strategy preserves
the pretrained backstepping knowledge while allowing adaptation to the
nonlinear dynamics.

Three DeepONet instances are initialized from the same pretrained
parameters: one for the actor, one shared by the twin online critics, and
one shared by the twin target critics. The actor and online critic
DeepONet parameters are fine-tuned independently, while the target
critic DeepONet is updated through Polyak averaging.

% The SAC actor generates the stochastic control component according to
% \begin{equation}
% U_{\rm SAC}\sim\pi_{\phi}(\cdot|z_\pi),
% \end{equation}
% and the boundary input applied to the plant is the sum of the SAC actor output and the prior control, as defined in \eqref{eq:final_U}.
 
% Thus, the critic evaluates the composite boundary input $U_L$, rather than
% $U_{\rm SAC}$ alone.

% A transfer-learning strategy is further introduced for the DeepONet
% parameters. For each online DeepONet network
% $D\in\{\pi,Q\}$, its parameters are partitioned as
% \begin{equation}
% \theta_D=(\theta_{f,D},\theta_{s,D}),
% \label{eq:parameter_partition}
% \end{equation}
% where $\theta_{f,D}$ is frozen during SAC training and $\theta_{s,D}$ is
% fine-tuned using a reduced learning rate $\eta_s$. Specifically, the frozen
% part contains the fully connected layers except the last layer of the trunk
% network and the convolutional layers of the branch network, while the
% remaining layers are allowed to adapt. 

%Three DeepONet instances are initialized from the same pretrained parameters: one for the actor, one for the online critics, and one for the target critics. The actor and online critic networks are adapted independently, whereas the target networks are updated through Polyak averaging.

The training procedure is summarized in Algorithm~\ref{alg:tdpsac}.
\begin{algorithm}[t]
\caption{Training of Transfer DeepONet Prior SAC (TDPSAC)}
\label{alg:tdpsac}
\begin{algorithmic}[1]

\State Pretrain DeepONet to approximate the backstepping feedback
operator and obtain pretrained parameters $\theta_D$.

\State Initialize actor, critic, and target critic DeepONets using
$\theta_D$; freeze $\theta_{f,\pi}$ and $\theta_{f,Q}$.

\State Initialize SAC networks and replay buffer $\mathcal D$.

\For{each episode}

    \State Observe initial state $s^0$.

    \For{$k=0,1,\ldots,K-1$}

        \State Compute $z_\pi^k$ and $U_{\rm prior}^k$ using the
        actor-side DeepONet.

        \State Sample
        $U_{\rm SAC}^k\sim\pi_\phi(\cdot|z_\pi^k)$ and construct 
        \[
        U_L^k=\rho^kU_{\rm SAC}^k+(1-\rho^k)U_{\rm prior}^k .
        \]

        \State Apply $U_L^k$, collect transition, and store it in
        $\mathcal D$.

        \State Update SAC actor and twin critics, together with
        $\theta_{s,\pi}$ and $\theta_{s,Q}$.

        \State Update target critics by Polyak averaging.

    \EndFor

\EndFor

\end{algorithmic}
\end{algorithm}

The closed-loop transition of the discretized Saint-Venant system is
written as
\begin{equation}
s^{k+1}=\mathcal F_{\Delta t}(s^k,U_L^k),
\label{eq:discrete_transition}
\end{equation}
where $\mathcal F_{\Delta t}$ denotes the numerical flow map over one
control interval $\Delta t$.

% Let $s^k$ denote the spatially discretized Saint-Venant state at the $k$th
% control instant. The closed-loop transition induced by the composite control
% can be written compactly as
% \begin{equation}
% s^{k+1}
% =
% \mathcal{F}_{\Delta t}(s^k,U_L^k),
% \label{eq:discrete_transition}
% \end{equation}
% where $\mathcal{F}_{\Delta t}$ denotes the numerical evolution of the
% Saint-Venant system \eqref{eq:fx}-\eqref{eq:H_bnd} over one control interval $\Delta t$.

The reward is designed to penalize the state tracking errors, temporal variations of the
states, deviations of the control input from its steady-state value, and
control variations,
% \begin{equation}
% \begin{aligned}
% r^k=-\Bigg[
% &c_1\left\|\frac{H^k-\bar H}{H_{\mathrm{scale}}}\right\|_N^2
% +c_2\left\|\frac{Q^k-\bar Q}{Q_{\mathrm{scale}}}\right\|_N^2
% \\
% &+c_3\left\|\frac{H^{k+1}-H^k}{H_{\mathrm{scale}}}\right\|_N^2
% +c_4\left\|\frac{Q^{k+1}-Q^k}{Q_{\mathrm{scale}}}\right\|_N^2
% \\
% &+c_5\left|\frac{U_L^k-\bar U_L}{U_{\mathrm{scale}}}\right|^2
% +c_6\left|\frac{U_L^k-U_L^{k-1}}
% {U_{\mathrm{scale}}}\right|^2
% \Bigg],
% \label{eq:reward_tdpsac}
% \end{aligned}
% \end{equation}
\begin{equation}
\begin{aligned}
r^k=-\Bigg[
&c_1\left(
\left\|H^k-\bar H\right\|_N^2
+\left\|Q^k-\bar Q\right\|_N^2
\right) \\
&+c_2\left(
\left\|H^{k+1}-H^k\right\|_N^2
+\left\|Q^{k+1}-Q^k\right\|_N^2
\right) \\
&+c_3\left|
\frac{U_L^k-\bar U_L}{U_{\mathrm{scale}}}
\right|^2
+c_4\left|
\frac{U_L^k-U_L^{k-1}}{U_{\mathrm{scale}}}
\right|^2
\Bigg].
\label{eq:reward_tdpsac}
\end{aligned}
\end{equation}
where $\bar H$, $\bar Q$, and $\bar U_L$ denote the desired steady-state
profiles and desired valve position at $x=L$, respectively. 
$c_i>0$, $i=1,\ldots4$, are weighting coefficients.
 The normalized spatial norm is defined, for example
as
\begin{equation}
\|H^k-\bar H\|_N^2
=\Delta x\sum_{m=1}^{N}\frac{1}{H_{\rm{scale}}}
(H_m^k-\bar H_m)^2,
\label{eq:discrete_norm}
\end{equation}
with the same definition used for $Q$. The quantities $H_{\mathrm{scale}}$, $Q_{\mathrm{scale}}$, and
$U_{\mathrm{scale}}$ are normalization factors for the corresponding
variables.

The SAC updates follow the standard off-policy actor--critic procedure, except
that the critic evaluates the composite action $U_L$ in
\eqref{eq:final_U} rather than the SAC action alone.

\section{SIMULATION}
% \begin{table*}[!t]
% \centering
% \caption{Controllers used for comparison and ablation}
% \label{tab:controller_defs}
% \scriptsize
% \setlength{\tabcolsep}{2.5pt}
% \renewcommand{\arraystretch}{1.15}
% \begin{tabular*}{\textwidth}{@{\extracolsep{\fill}}p{0.14\textwidth}p{0.18\textwidth}p{0.21\textwidth}p{0.25\textwidth}p{0.14\textwidth}@{}}
% \hline
% Controller 
% & DeepONet as feature encoder 
% & DeepONet as control prior 
% & Trainable/frozen setting 
% & Ablation purpose \\
% \hline

% BS
% & No
% & Analytical backstepping law
% & No learning
% & Classical baseline \\

% SAC
% & No
% & No
% & No transfer
% & Pure SAC baseline \\

% SAC+BC
% & No
% & Independent DeepONet action prior
% & Independent DeepONet action prior fully frozen
% & Effect of frozen action prior only \\

% SAC-NOE
% & Actor/Critic DeepONet feature encoders
% & No
% & All DeepONet feature encoder trainable
% & Effect of DeepONet features only \\

% SAC-NOE-S
% & Actor/Critic DeepONet encoders
% & Actor DeepONet control branch
% & All DeepONet parameters trainable
% & Effect of fully trainable DeepONet control \\

% TSAC-FNO
% & Actor/Critic DeepONet feature encoders
% & Independent DeepONet action prior
% & Actor/Critic DeepONets partially frozen; independent DeepONet action prior fully frozen
% & Effect of additional frozen independent prior \\

% TSAC
% & Actor/Critic DeepONet encoders
% & Actor DeepONet control branch
% & Actor/Critic DeepONets partially frozen and fine-tuned
% & Proposed transfer strategy \\
% \hline
% \end{tabular*}
% \end{table*}

The simulations are conducted on the Saint-Venant canal model using the Sambre River parameters in~\cite{dosSantos2008boundary}, as listed in Table~\ref{tab:saint_venant_params}. The numerical solution of \eqref{eq:fx}--\eqref{eq:H_bnd} is obtained using the finite difference method, with spatial step $\Delta x=111.28~\mathrm{m}$, time step $\Delta t=5~\mathrm{s}$, and a simulation horizon of $T=666.67~\mathrm{min}$. The source code is available on GitHub (\url{https://github.com/wangchenchen-gif/Saint-Venant-DeepONet-RL}).

%The simulations are conducted on the Saint-Venant canal model using the Sambre River parameters in~\cite{dosSantos2008boundary}, as listed in Table~\ref{tab:saint_venant_params}. The numerical solution of \eqref{eq:fx}--\eqref{eq:H_bnd} is obtained using
%the finite difference method, by setting spatial step $\Delta x=111.28~\mathrm{m}$, time step $\Delta t=5~\mathrm{s}$, and a simulation horizon of $T=666.67~\mathrm{min}$. The source code is available on GitHub (\url{https://github.com/wangchenchen-gif/Saint-Venant-DeepONet-RL}).
\begin{table}[h]
	\centering
	\caption{Simulation Parameters for the Saint-Venant Model}
	\label{tab:saint_venant_params}
	\scriptsize
	\setlength{\tabcolsep}{3.0pt}
	\renewcommand{\arraystretch}{1.12}
	\begin{tabular}{p{0.16\columnwidth}p{0.26\columnwidth}p{0.46\columnwidth}}
		\hline
		Parameter & Value & Meaning\\
		\hline
		$L$ & $11239~\mathrm{m}$ & Channel length\\
		$B$ & $40~\mathrm{m}$ & Channel width\\
		$g$ & $9.81~\mathrm{m/s^2}$ & Gravity\\
		$I$ & $7.92\times 10^{-5}$ & Bottom slope\\
		$K$ & $33~\mathrm{m^{1/3}s^{-1}}$ & Strickler coefficient\\
		$n_M$ & $1/33~\mathrm{s~m^{-1/3}}$ & Manning coefficient\\
		$\mu_0$ & $0.4$ & Upstream discharge coefficient\\
		$\mu_L$ & $0.4$ & Downstream discharge coefficient\\
		$z_{\rm up}$ & $5.00~\mathrm{m}$ & Upstream reservoir level\\
		$\bar U_0$ & $0.1548~\mathrm{m}$ & Equilibrium upstream gate opening\\
		$\bar U_L$ & $0.30~\mathrm{m}$ & Equilibrium downstream gate opening\\
		$h_s$ & $4.094~\mathrm{m}$ & Downstream sill height\\
		\hline
	\end{tabular}
\end{table}

For DeepONet pretraining, the initial water-depth and discharge profiles are generated by truncated Chebyshev expansions,
\begin{equation}
\begin{aligned}
H_{\ell}^{0}(x)
&=\bar H(x)+\sum_{j=0}^{M} a_{\ell,j}
T_j\left(\frac{2x}{L}-1\right),\\
Q_{\ell}^{0}(x)
&=\bar Q(x)+\sum_{j=0}^{M} b_{\ell,j}
T_j\left(\frac{2x}{L}-1\right),
\end{aligned}
\label{eq:random_initial_profiles}
\end{equation}
where $T_j$ is the $j$th Chebyshev polynomial, and $a_{\ell,j}$ and $b_{\ell,j}$ are independently sampled from zero-mean Gaussian distributions with standard deviations proportional to $(j+1)^{-2}$. We set $M=6$ and limit the maximum perturbation amplitudes to $0.30~\mathrm{m}$ for water depth and $2.50~\mathrm{m^3/s}$ for discharge. Only trajectories converging under the Backstepping controller are retained, resulting in $10^5$ training samples.

During RL training, the truncated-Chebyshev initial conditions are resampled at each episode, and each controller is trained for $10^5$ environment interactions. The replay-buffer and mini-batch sizes are $10^5$ and $256$, respectively, with SAC learning rate $10^{-4}$, DeepONet fine-tuning rate $\eta_s=10^{-6}$, Polyak coefficient $\tau=0.005$, and discount factor $\gamma=0.995$. The blending factor defined in \eqref{eq:final_U} is set to
$\rho=0.1$, $0.5$, and $0.9$ over $0$--$10\%$, $10$--$20\%$,
and $20$--$100\%$ of training, respectively, providing a staged
transition from prior-guided control to SAC-dominated adaptation. The reward weights are $c_1=1.00$, $c_2=0.10$, $c_3=0.05$, and $c_4=0.10$, with $H_{\mathrm{scale}}=0.30~\mathrm{m}$, $Q_{\mathrm{scale}}=3.00~\mathrm{m^3/s}$, and $U_{\mathrm{scale}}=0.10~\mathrm{m}$. All experiments were conducted on a workstation equipped with an Intel Core i9-13900KF CPU and an NVIDIA GeForce RTX 4090 GPU.
 
\subsection{Backstepping Prior Learning}

The analytical Backstepping controller is first used to generate control labels for the sampled initial profiles in \eqref{eq:random_initial_profiles}. The resulting state--control pairs are then used to pretrain the DeepONet so that its prior-control output $U_{\mathrm{prior}}$ approximates the Backstepping control law $U_{\mathrm{bs}}$.

% To evaluate the approximation accuracy, we consider the initial condition
% [
% H(0,x)=3.75+\frac{4.65-3.75}{L}x~\mathrm{m},\qquad
% Q(0,x)=10~\mathrm{m^3/s}.
% ]
The pretrained DeepONet achieves an average $L^2$ norm approximation error of less than $10^{-4}$ over 50 evaluations,  
% $
% \|U_{\mathrm{prior}}-U_{\mathrm{bs}}\|_{L^2}<10^{-4},
% $
showing that the backstepping control mapping is accurately captured and can serve as a reliable prior for subsequent reinforcement learning.

\subsection{Learning Efficiency and Ablation Study}

\begin{table*}[!t]
\centering
\caption{Controllers used for comparison and ablation studies.}
\label{tab:controller_defs}
\scriptsize
\setlength{\tabcolsep}{2.5pt}
\renewcommand{\arraystretch}{1.15}
\begin{tabular*}{\textwidth}{
@{\extracolsep{\fill}}
>{\raggedright\arraybackslash}p{0.09\textwidth}
>{\raggedright\arraybackslash}p{0.20\textwidth}
>{\raggedright\arraybackslash}p{0.20\textwidth}
>{\raggedright\arraybackslash}p{0.43\textwidth}
% >{\raggedright\arraybackslash}p{0.15\textwidth}
@{}
}
\hline
\textbf{Controller}
& \textbf{Feature Extractor}
& \textbf{Prior Control $U_{\mathrm{prior}}$}
& \textbf{RL Training Strategy} \\
\hline

BS
& None
& Analytical Backstepping
& None \\

SAC
& None
& None 
& Training from scratch \\

SAC-FPC
& None
& Pretrained DeepONet
& Prior DeepONet frozen \\

DO-SAC
& Pretrained DeepONet
& None
& DeepONet fully trainable\\

DPSAC
& Pretrained DeepONet
& Pretrained DeepONet
& Actor extractor and prior share the same DeepONet with fully trainable \\

TDPSAC-FPC
& Pretrained DeepONet
& Pretrained DeepONet
& Extractor transfer training; prior frozen\\

TDPSAC
& Pretrained DeepONet 
& Pretrained DeepONet 
& Actor extractor and prior share the same DeepONet with transfer training \\

\hline
\end{tabular*}
\end{table*}
To evaluate the contribution of each component in TDPSAC, six learning-based controllers are compared, as summarized in Table~\ref{tab:controller_defs}. These baselines separately examine the effects of DeepONet feature encoding, prior-control guidance, transfer learning, and prior adaptation.

\begin{figure}[htbp]
    \centering
    \includegraphics[width=0.7\linewidth]{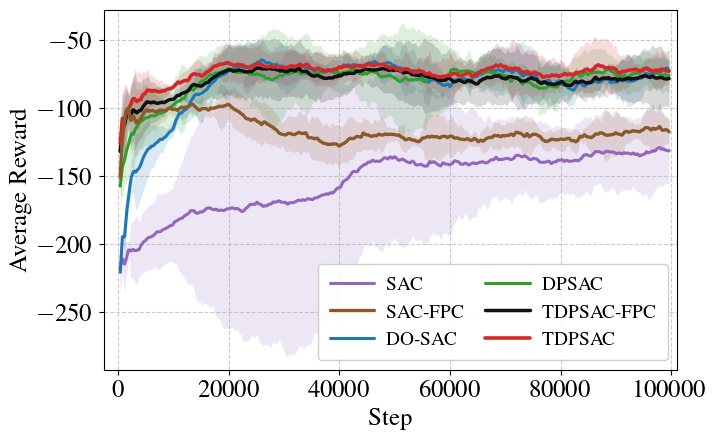}
    \caption{RL training reward curves. Solid lines and shaded regions denote the mean and 95\% confidence interval over five seeds.}
    \label{fig:re}
\end{figure}

%Fig.~\ref{fig:re} compares the mean episodic returns over five random seeds. Since the reward is the negative control cost, returns closer to zero indicate better performance. TDPSAC, TDPSAC-FPC, and DPSAC improve rapidly at the early stage, showing that the pretrained DeepONet representation and control prior reduce ineffective exploration. DO-SAC also benefits from DeepONet features but converges more slowly without prior-action guidance, while SAC performs worst because both the representation and policy are learned from scratch. SAC-FPC performs well initially due to the pretrained prior, but its frozen prior cannot adapt to the evolving policy and state distribution, causing increasing mismatch with the SAC correction and degraded later-stage performance.
%Among all learning-based controllers, TDPSAC achieves the highest and most stable return. Its advantage over DPSAC shows the benefit of partial fine-tuning, while its advantage over TDPSAC-FPC shows that adapting the prior-control branch during RL further improves performance.
Fig.~\ref{fig:re} compares the mean episodic returns over five random seeds. Since the reward is defined as the negative control cost, a return closer to zero indicates better performance. TDPSAC, TDPSAC-FPC, and DPSAC improve rapidly during the early training stage, showing that the pretrained DeepONet representation and prior control reduce ineffective exploration. DO-SAC also benefits from DeepONet feature extraction, but converges more slowly without prior-action guidance, while SAC performs worst because both the state representation and control policy are learned from scratch. SAC-FPC performs well initially because the pretrained prior provides effective action guidance. However, since the prior remains frozen during RL training, it cannot adapt to the evolving policy and state distribution, which may lead to increasing mismatch between the prior action and the SAC correction and hence degraded performance at later stages.

Among all learning-based controllers, TDPSAC achieves the highest and most stable return. The comparison with DPSAC indicates the benefit of retaining transferred knowledge through partial fine-tuning, while the comparison with TDPSAC-FPC shows that allowing the prior-control branch to adapt during RL further improves performance. 

\subsection{Closed-Loop Regulation Performance}
%To avoid abrupt variations in the downstream gate command, a first-order
%low-pass filter is applied to the control input \cite{Itahashi2024}, i.e.,
%$U_{L,f}^{k}=(1-\alpha_f)U_{L,f}^{k-1}+\alpha_f U_L^{k}$,
%where $\alpha_f\in(0,1]$ is the filtering coefficient. The filtered
%signal $U_{L,f}^{k}$ is applied to the downstream gate, with
%$\alpha_f=0.7$ in the simulations.

\begin{figure}[htbp]
	\centering
	\subfloat[]{\includegraphics[width=0.24\textwidth]{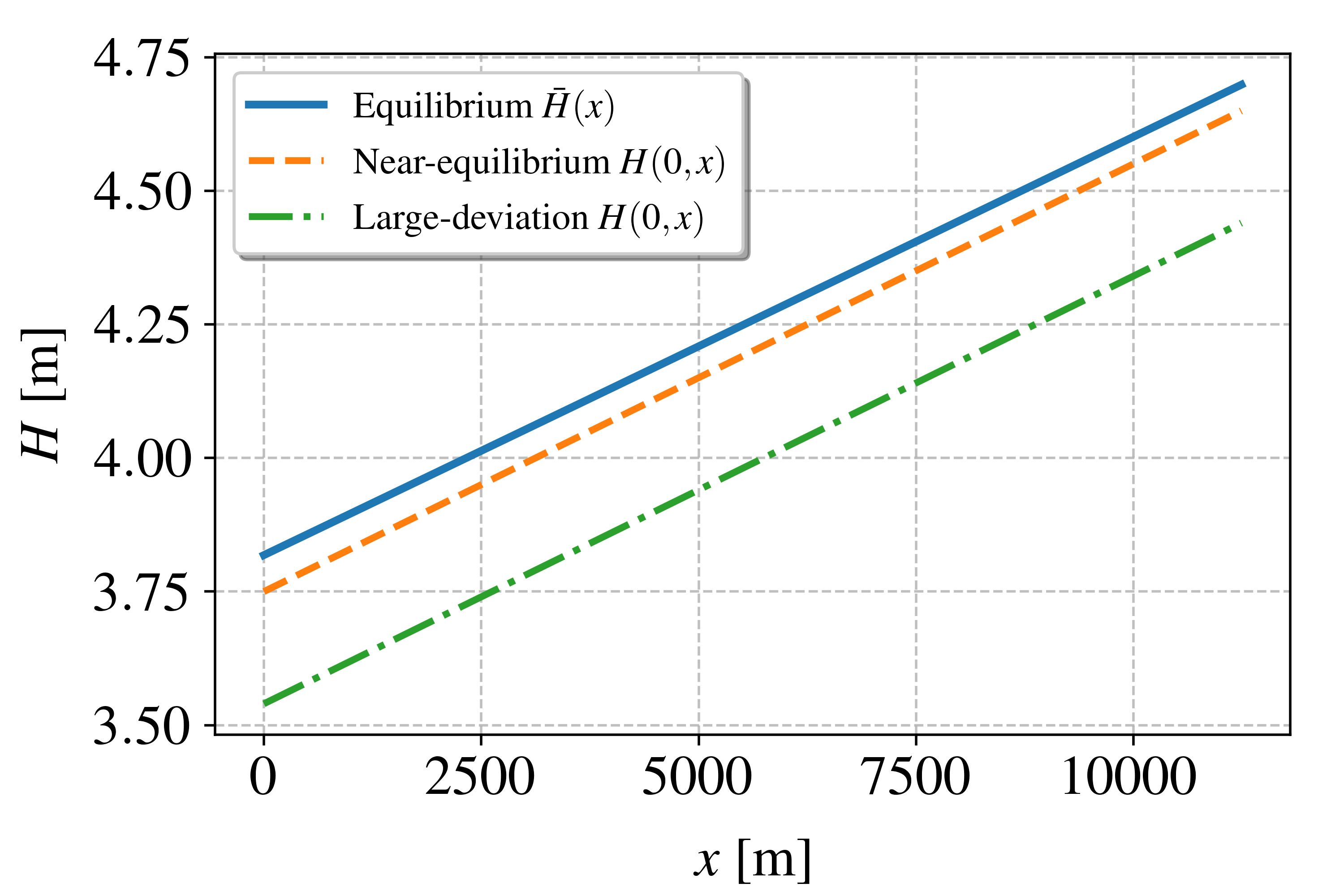}}
	\subfloat[]{\includegraphics[width=0.24\textwidth]{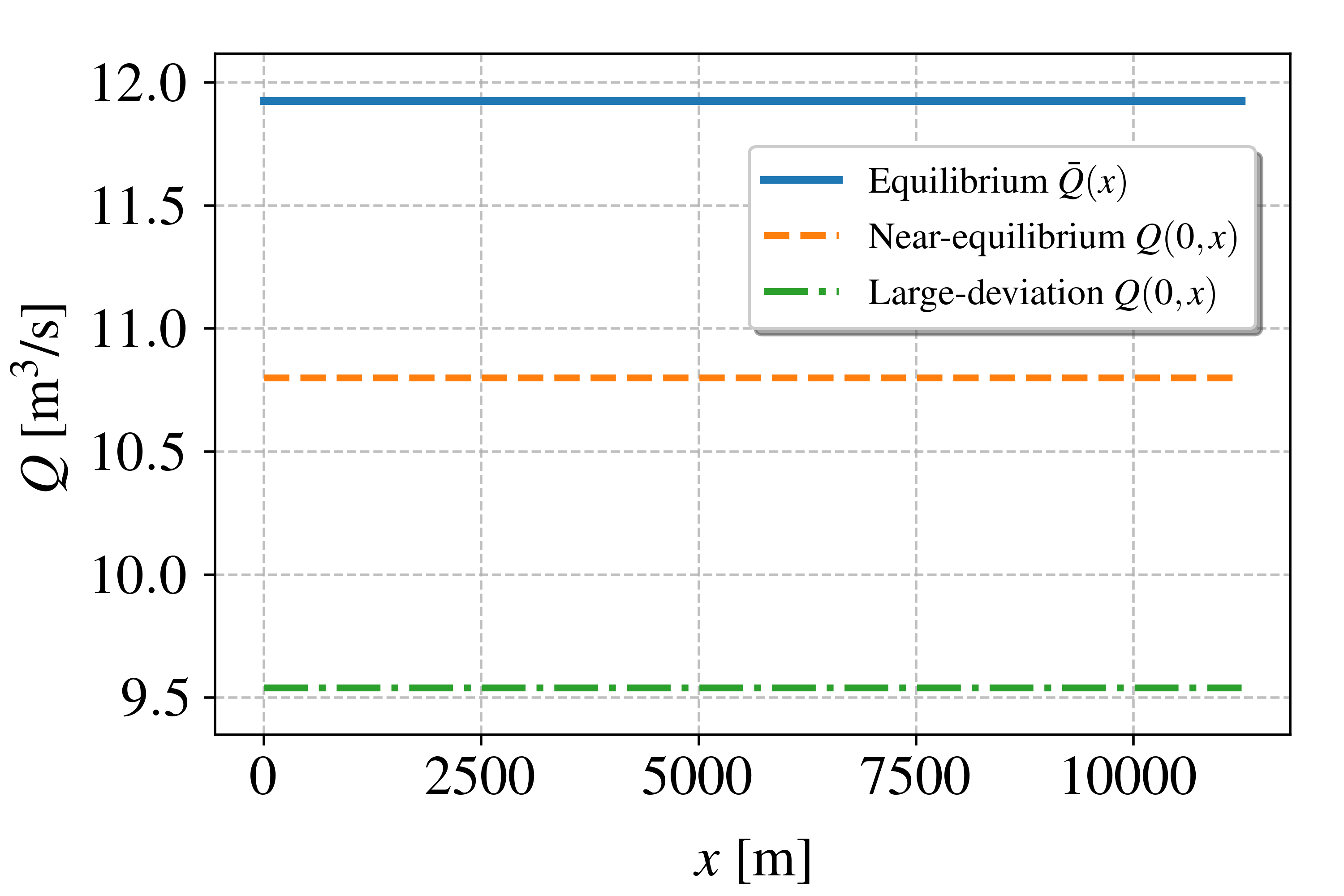}}
	\caption{ 
		Equilibrium of the Saint-Venant canal system obtained by solving
the steady-state Saint-Venant equation.}
	\label{fig:equilibrium_profiles}
\end{figure}

\begin{figure}[htbp]
\centering
\subfloat[]{\includegraphics[width=0.24\textwidth]{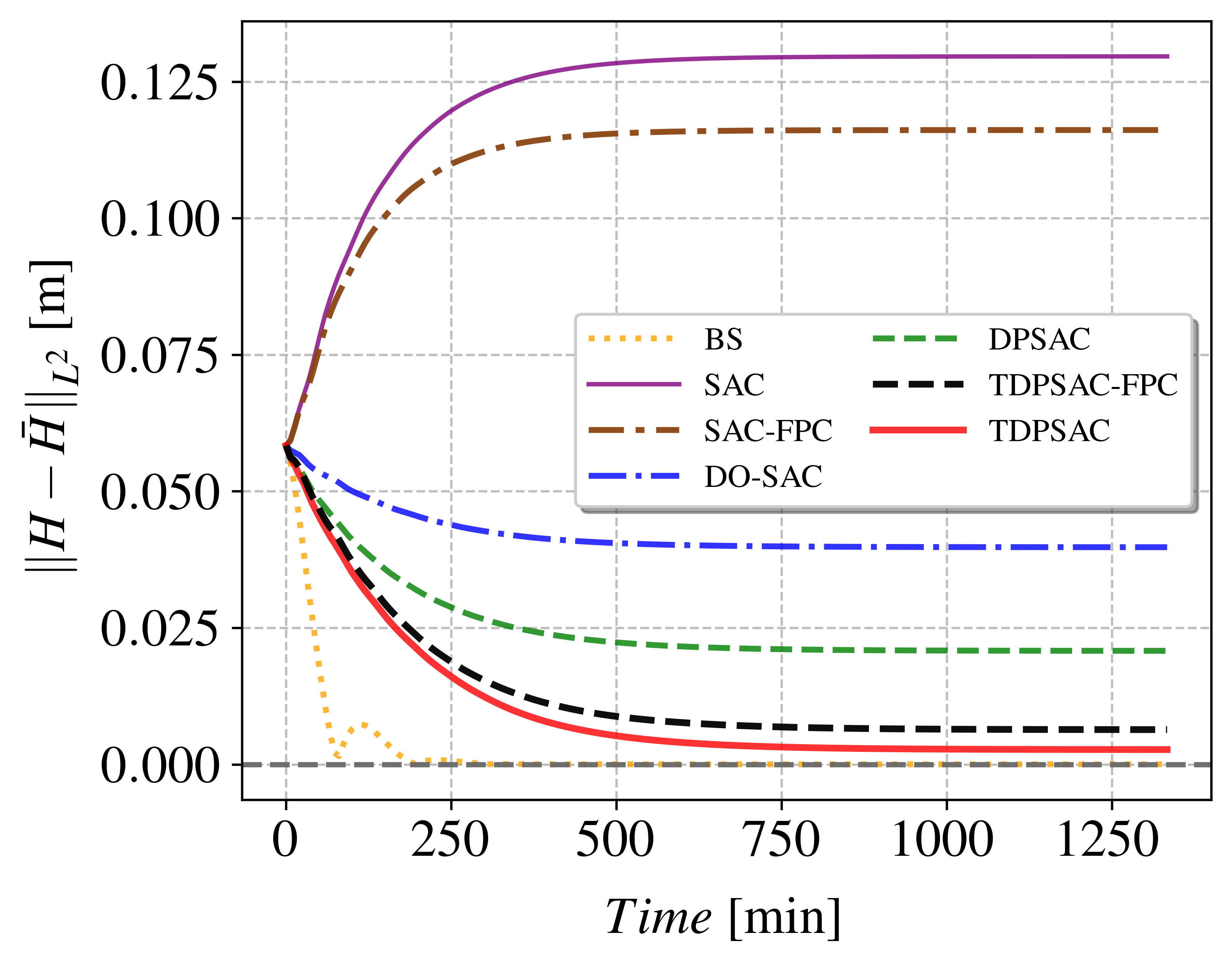}}
\subfloat[]{\includegraphics[width=0.24\textwidth]{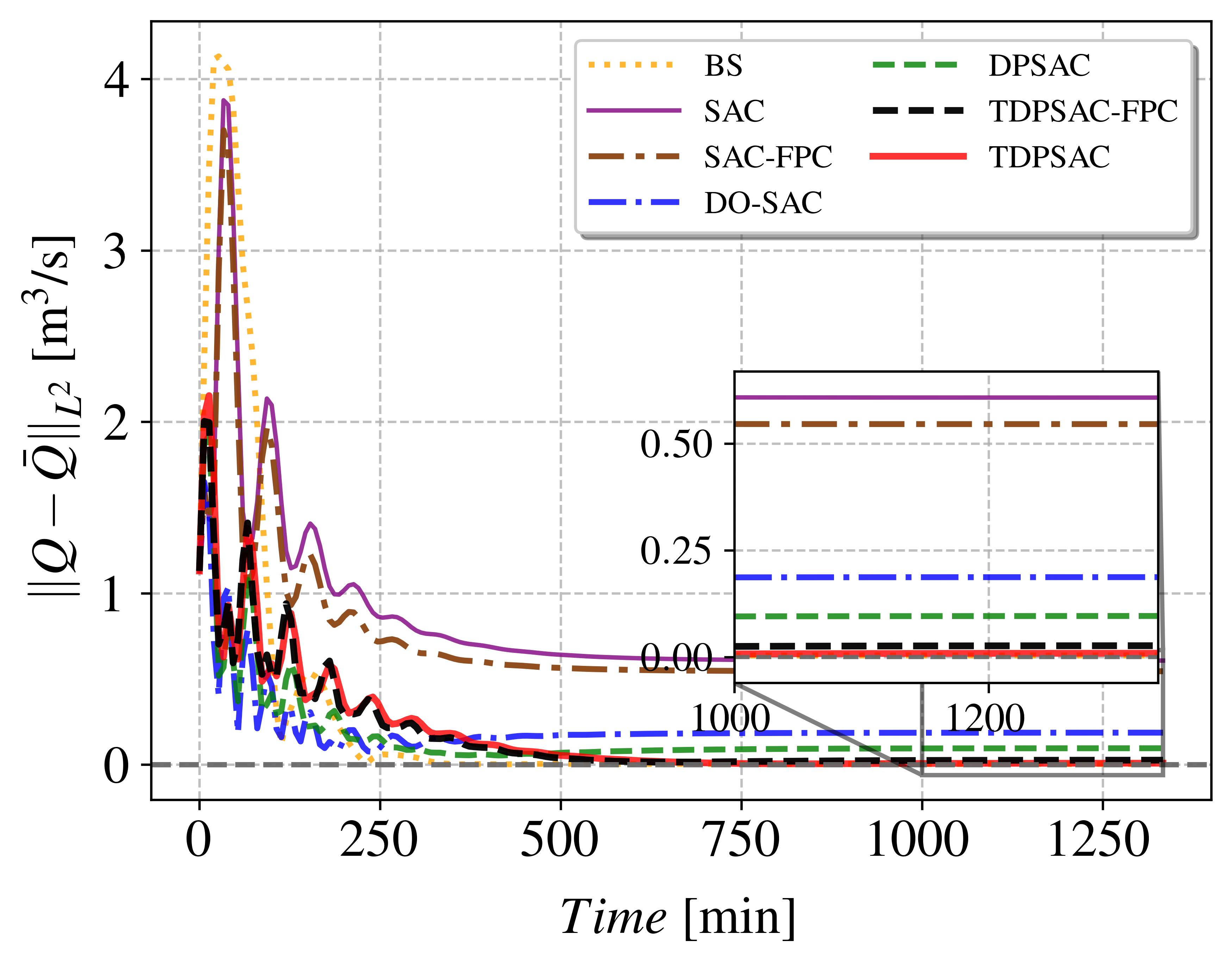}}
\caption{State $L^2$-norm errors of the compared controllers in the near-equilibrium case:
(a) $\|H-\bar H\|_{L^2}$ and (b) $\|Q-\bar Q\|_{L^2}$,
with $H(0,x)=3.75+(4.65-3.75)x/L~\mathrm{m}$ and
$Q(0,x)=10.8~\mathrm{m^3/s}$.}
\label{fig:L}
\end{figure}

A first-order low-pass filter with coefficient $0.7$ is applied to the downstream gate command to suppress abrupt control variations \cite{Itahashi2024}.
We first evaluate the controllers under near-equilibrium operating conditions. The relative initial deviations are defined as
$\epsilon_H=\|H(0,\cdot)-\bar H(\cdot)\|_{L^2}/\|\bar H(\cdot)\|_{L^2}\times100\%$
and
$\epsilon_Q=\|Q(0,\cdot)-\bar Q(\cdot)\|_{L^2}/\|\bar Q(\cdot)\|_{L^2}\times100\%$. Since the desired equilibrium for state $H$ is approximately linear, as shown in Fig.~\ref{fig:equilibrium_profiles}, we chose a linear initial water-depth profile with two ends 
$H(0,0)=3.75~\mathrm{m}$ and 
$H(0,L)=4.65~\mathrm{m}$, together with 
$Q(0,x)=10.8~\mathrm{m^3/s}$. These initial conditions correspond to 
$\epsilon_H=1.4\%$ and $\epsilon_Q=9.5\%$.
As shown in Fig.~\ref{fig:L}, all controllers successfully regulate the system under this near-equilibrium initial condition. The analytical Backstepping controller achieves the smallest steady-state errors, as expected from its model-based design around the nominal equilibrium. However, it exhibits relatively larger transient overshoots. In contrast, the learning-based controllers are trained using a multi-objective reward that balances tracking accuracy, temporal smoothness, and control effort, and therefore exhibit smaller overshoots. Among the learning-based methods, TDPSAC achieves the best overall performance, with the smallest steady-state errors and overshoot. TDPSAC-FPC ranks second, while DPSAC and DO-SAC outperform SAC and SAC-FPC, demonstrating the benefits of transferred DeepONet features and prior-control guidance.

\begin{figure}[htbp]
\centering
\subfloat[]{\includegraphics[width=0.24\textwidth]{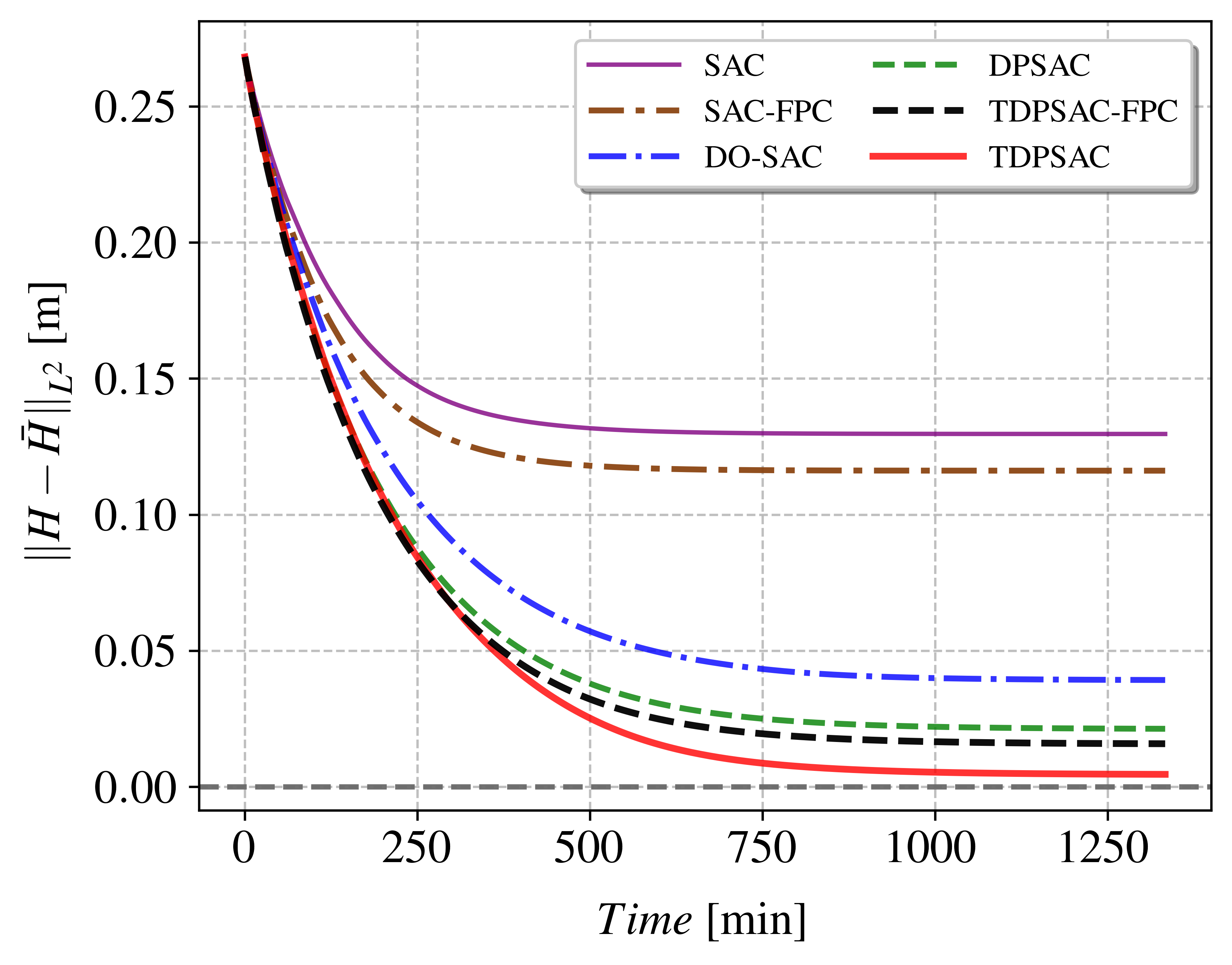}}
\subfloat[]{\includegraphics[width=0.24\textwidth]{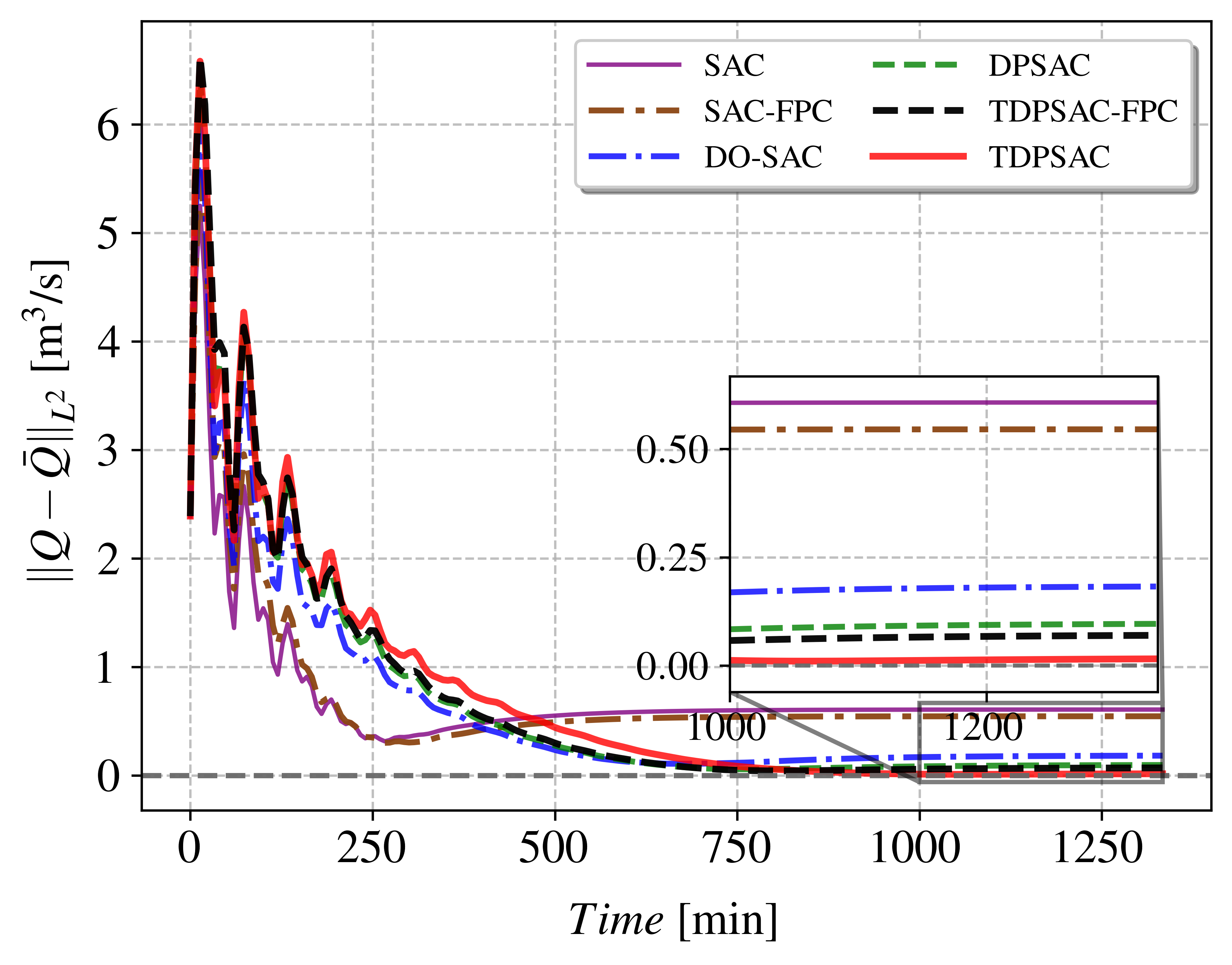}}
\caption{State $L^2$-norm errors of the compared controllers in the large-deviation case:
(a) $\|H-\bar H\|_{L^2}$ and (b) $\|Q-\bar Q\|_{L^2}$,
with $H(0,x)=3.54+(4.44-3.54)x/L~\mathrm{m}$ and
$Q(0,x)=9.54~\mathrm{m^3/s}$.}
\label{fig:Lbig}
\end{figure}

We next consider a larger initial deviation with
$H(0,0)=3.54~\mathrm{m}$,
$H(0,L)=4.44~\mathrm{m}$, and
$Q(0,x)=9.54~\mathrm{m^3/s}$, corresponding to
$\epsilon_H=6.3\%$ and $\epsilon_Q=20\%$.
In this case, the analytical Backstepping controller fails to maintain closed-loop convergence, while all learning-based controllers remain stable. As shown in Fig.~\ref{fig:Lbig}, TDPSAC achieves the best performance among the learning-based controllers. The performance gaps mainly arise from differences in how pretrained backstepping knowledge is incorporated and preserved during reinforcement learning. TDPSAC-FPC ranks second since its frozen prior limits adaptation to the changing state distribution during RL. DPSAC updates all DeepONet parameters and may therefore weaken the retention of pretrained control knowledge, whereas DO-SAC uses DeepONet only for feature extraction without explicit prior-action guidance. SAC lacks pretrained control knowledge, while SAC-FPC relies on a fixed prior that cannot adapt during training.

To statistically assess the regulation performance, 50 trials are conducted for each case, with the initial conditions randomly sampled within the prescribed relative deviation ranges $\epsilon_H$ and $\epsilon_Q$, respectively. For both $H$ and $Q$, the samples are distributed on both sides of the corresponding equilibrium profiles, with 25 samples above and 25 below the equilibrium. The results are summarized in Table~\ref{tab:performance_different_initial_conditions}. Under near-equilibrium initial conditions, all controllers successfully stabilize the nonlinear system. The Backstepping controller achieves the smallest steady-state errors, while the learning-based controllers show smaller transient overshoots. 
%Among them, TDPSAC achieves the lowest steady state error $\epsilon_H$, $\epsilon_Q$, and $M_Q$. %Among the learning-based controllers, TDPSAC attains the lowest $\epsilon_H$ and $\epsilon_Q$, together with a relatively small discharge overshoot $M_Q$.

Under large initial deviations, the Backstepping controller cannot guarantee closed-loop convergence, whereas all learning-based controllers successfully stabilize the system. In particular, TDPSAC achieves the lowest steady-state errors among all learning-based controllers. These results indicate that the transferred DeepONet representation and adaptive prior control preserve useful knowledge from the Backstepping controller, while SAC adaptation further extends the effective regulation range beyond that of the conventional analytical controller. 

The RL training times of the learning-based controllers are comparable, ranging from 32 to 41 min. TDPSAC requires only a modest additional training cost while achieving substantially improved regulation performance.

%Overall, the advantage of TDPSAC does not arise from either the pretrained representation or the prior-control branch alone, but from their joint use with partial parameter fine-tuning and online SAC policy correction.
\begin{table}[h]
\centering
\caption{Closed-loop performance under different initial conditions.}
\label{tab:performance_different_initial_conditions}
\setlength{\tabcolsep}{3.0pt}
\renewcommand{\arraystretch}{1.12}
\resizebox{\columnwidth}{!}{%
\begin{tabular}{
@{}
l
c
c
c
c
c
@{}
}
\toprule
Controller
& \makecell{$\epsilon_H$\\$[\times10^{-3}]$}
& \makecell{$\epsilon_Q$\\$[\times10^{-3}]$}
& \makecell{$M_H$\\$[\%]$}
& \makecell{$M_Q$\\$[\%]$}
& \makecell{$t_{\mathrm{RL}}$\\$[\mathrm{min}]$}
\\
\midrule

\multicolumn{6}{@{}l@{}}{%
\textbf{Case I: Near-equilibrium initial conditions
$\left(\epsilon_H\in[1\%,2\%],\,\epsilon_Q\in[9\%,10\%]\right)$}}
\\
\midrule

BS
& $\mathbf{0.0181 \pm 6.06\times 10^{-6}}$
& $\mathbf{4.14 \pm 4.61\times 10^{-5}}$
& $1.09 \pm 0.842$
& $42.1 \pm 25.5$
& --
\\

SAC
& $130 \pm 0.00508$
& $608 \pm 0.138$
& $\mathbf{0.943 \pm 0.965}$
& $60.9 \pm 26.6$
& 32
\\

SAC-FPC
& $116 \pm 0.00195$
& $545 \pm 0.0584$
& $\mathbf{0.943 \pm 0.965}$
& $59.0 \pm 28.1$
& 37
\\

DO-SAC
& $39.7 \pm 0.0683$
& $188 \pm 1.37$
& $\mathbf{0.943 \pm 0.965}$
& $35.3 \pm 17.4$
& 36
\\

DPSAC
& $20.8 \pm 0.0464$
& $96.5 \pm 0.915$
& $0.945 \pm 0.967$
& $30.5 \pm 13.8$
& 39
\\

TDPSAC-FPC
& $6.35 \pm 0.0455$
& $26.8 \pm 1.02$
& $0.946 \pm 0.968$
& $28.5 \pm 12.0$
& 40
\\

TDPSAC
& $2.72 \pm 0.0261$
& $8.78 \pm 0.614$
& $0.947 \pm 0.969$
& $\mathbf{28.0 \pm 11.6}$
& 41
\\

\midrule
\multicolumn{6}{@{}l@{}}{%
\textbf{Case II: Large-deviation initial conditions
$\left(\epsilon_H\in[6\%,7\%],\,\epsilon_Q\in[20\%,21\%]\right)$}}
\\
\midrule

SAC
& $130 \pm 0.0200$
& $608 \pm 0.545$
& $\mathbf{3.77 \pm 3.81}$
& $106 \pm 77.7$
& 32
\\

SAC-FPC
& $116 \pm 0.00849$
& $545 \pm 0.255$
& $\mathbf{3.77 \pm 3.81}$
& $105 \pm 79.1$
& 37
\\

DO-SAC
& $39.1 \pm 0.349$
& $185 \pm 7.00$
& $\mathbf{3.77 \pm 3.81}$
& $102 \pm 60.9$
& 36
\\

DPSAC
& $21.2 \pm 0.382$
& $98.8 \pm 6.80$
& $\mathbf{3.77 \pm 3.81}$
& $98.1 \pm 58.1$
& 39
\\

TDPSAC-FPC
& $15.7 \pm 0.355$
& $71.8 \pm 6.21$
& $\mathbf{3.77 \pm 3.81}$
& $\mathbf{95.4 \pm 58.3}$
& 40
\\

TDPSAC
& $\mathbf{4.57 \pm 0.276}$
& $\mathbf{16.9 \pm 4.35}$
& $\mathbf{3.77 \pm 3.81}$
& $98.4 \pm 56.3$
& 41
\\

\bottomrule
\end{tabular}%
}

\vspace{2pt}
\begin{minipage}{\columnwidth}
\scriptsize
\textit{Note:}
Values are reported as mean $\pm$ standard deviation across 50 trials.
$\epsilon_H$ and $\epsilon_Q$ denote the steady-state $L^2$-norm errors of the water-depth and discharge profiles, respectively.
$M_H$ and $M_Q$ denote the maximum positive relative overshoots over space and time; and $t_{\mathrm{RL}}$ denotes the RL training time.
The backstepping controller is excluded from the large-deviation case because it fails to converge.
\end{minipage}
\end{table}

\section{Conclusion}
\label{sec:conclusion}
This letter presents a backstepping-informed reinforcement learning framework for regulating the nonlinear Saint-Venant canal system. The numerical results show that the proposed TDPSAC achieves a favorable balance between local regulation accuracy, transient performance, and adaptability to initial conditions away from the equilibrium. Ablation studies further demonstrate the benefit of combining transferred backstepping knowledge with adaptive policy learning. 
These results indicate that embedding analytical control priors into RL can extend the effective operating range of model-based PDE controllers and improve learning efficiency. Future work will address robustness to persistent disturbances and slowly varying operating conditions.

\bibliographystyle{IEEEtran}
\bibliography{Ref}
\end{document}